\documentclass[twoside,11pt]{article}

\usepackage{graphicx, subfig}
\usepackage[top=1in,bottom=1in,left=1in,right=1in]{geometry}
\usepackage[sort&compress]{natbib} 
\usepackage{amsfonts, amsmath, amssymb, amsthm, constants, bbm}
\usepackage{MnSymbol}
\usepackage{enumitem}
\usepackage[mathscr]{eucal}
\usepackage{booktabs}

\usepackage{color}
\definecolor{darkred}{RGB}{100,0,0}
\definecolor{darkgreen}{RGB}{0,100,0}
\definecolor{darkblue}{RGB}{0,0,150}

\usepackage{hyperref}
\hypersetup{colorlinks=true, linkcolor=darkred, citecolor=darkgreen, urlcolor=darkblue}
\usepackage{url}

\newtheorem{thm}{Theorem}

\newtheorem{lem}{Lemma}

\theoremstyle{remark}
\newtheorem{Def}{Definition}

\theoremstyle{definition}
\newtheorem{con}{Contribution}

\def\beq{\begin{equation}} 
\def\eeq{\end{equation}}
\def\beqn{\begin{eqnarray*}}
\def\eeqn{\end{eqnarray*}}
\def\Bitem{\begin{itemize}\setlength{\itemsep}{.2in}}
\def\bitem{\begin{itemize}\setlength{\itemsep}{.05in}}
\def\eitem{\end{itemize}}
\def\Benum{\begin{enumerate}\setlength{\itemsep}{.2in}}
\def\benum{\begin{enumerate}\setlength{\itemsep}{.05in}}
\def\eenum{\end{enumerate}}
\def\bmult{\begin{multline*}}
\def\emult{\end{multline*}}
\def\bcenter{\begin{center}}
\def\ecenter{\end{center}}
\def\bframe{\begin{frame}}
\def\eframe{\end{frame}}

\newcommand{\thmref}[1]{Theorem~\ref{thm:#1}}

\newcommand{\lemref}[1]{Lemma~\ref{lem:#1}}
\newcommand{\secref}[1]{Section~\ref{sec:#1}}
\newcommand{\figref}[1]{Figure~\ref{fig:#1}}
\newcommand{\tabref}[1]{Table~\ref{tab:#1}}

\newcommand{\defref}[1]{Definition~\ref{def:#1}}

\def\cI{\mathcal{I}}
\def\cJ{\mathcal{J}}

\def\bS{\mathbf{S}}

\def\bX{\mathbf{X}}

\def\bx{\mathbf{x}}

\def\bbR{\mathbb{R}}

\newcommand{\E}{\operatorname{\mathbb{E}}}
\renewcommand{\P}{\operatorname{\mathbb{P}}}

\def\1{\mathbbm{1}}

\def\X{\bX}

\def\scan{\textsc{scan}}

\definecolor{purple}{rgb}{0.4,.1,.9}

\urldef\jgurl\url{jiaqig@amazon.com}
\urldef\ylurl\url{yuchao.liu@microsoft.com}

\begin{document}
	\thispagestyle{empty}
	
	\title{Distribution-Free, Size Adaptive Submatrix Detection with Acceleration}
	\author{
		Yuchao Liu\thanks{Microsoft Corporation --- \ylurl} 
		\and 
		Jiaqi Guo\thanks{Amazon Inc. \jgurl}
	}
	\date{}
	\maketitle

\begin{abstract}
Given a large matrix containing independent data entries, we consider the problem of detecting a submatrix inside the data matrix that contains larger-than-usual values. Different from previous literature, we do not have exact information about the dimension of the potential elevated submatrix. We propose a Bonferroni type testing procedure based on permutation tests, and show that our proposed test loses no first-order asymptotic power compared to tests with full knowledge of potential elevated submatrix. In order to speed up the calculation during the test, an approximation net is constructed and we show that Bonferroni type permutation test on the approximation net loses no power on the first order asymptotically.
\end{abstract}


\section{Introduction} \label{sec:intro}

Matrix type data are common in contemporary data analysis and have wide applications in biology, social sciences and other fields. In many situations, the row and column indexes represent individuals or units that could interact with each other, and the entries / data points evaluate the level of the iteration between row and column units. For example, in DNA chips analysis, the rows represent genes and columns represent situations. The entries inside the data matrix can be expression levels of some genes under some situations. See \citep{cheng2000biclustering} for a detailed introduction.

Bi-clustering, or co-clustering, or sometimes referred as simultaneous clustering, aims to find subsets of row and column indexes, such that entries inside the submatrix indexed by those units, are 'special'. This hidden structure is usually of interest to researchers and can be used for specific references. See \citep{tanay2005biclustering, charrad2011simultaneous, kriegel2009clustering} for surveys on this topic.

Before making inferences based on the bi-clustering result, one important question is that, does the result contain any information, or is it just a product of pure noise? This leads to the submatrix detection problem, which is a hypothesis testing procedure distinguishing data from pure noise and data containing some hidden structures.  

\subsection{Submatrix detection}

We consider the simplest case where there is one potential 'special' submatrix to discover. Furthermore, assume that entries are 'larger than usual' inside the submatrix. The observed data is denoted by $\bX = \{X_{ij}\}$ with $M$ rows and $N$ columns. The hypothesis testing problem is formulated as follows:
\begin{center}
$H_0$: $X_{ij}$ are IID for all $i\in [M]$ and $j\in [N]$ \\ 
vs\\
$H_1$: There exists $\cI\subset [M]$ and $\cJ \subset [N]$ such that $X_{ij}$ are stochastically larger when $(i,j) \in \cI \times \cJ$.
\end{center}
Here we make a notation that $[M] = \{1,2,\ldots,M\}$. To make the alternative hypothesis mathematically concrete, in $H_1$ we mean that for any $(i,j)\in \cI \times \cJ$ and $(k,l) \notin \cI \times \cJ$, and for any $c \in \mathbb{R}$, we have $\P(X_{ij} \geq c) \geq P(X_{kl} \geq c)$. And there exists $c = c_0\in \mathbb{R}$ that the equality does not hold.
 This problem is considered in \citep{butucea2013detection} with test statistics as 
\beq \label{sum}
\textsc{sum}(\X) = \sum_{i \in [M]} \sum_{j \in [N]} X_{ij},
\eeq
and 
\beq \label{scan}
\scan_{m,n}(\X)= \max_{\cI \subset [M], |\cI| = m} \quad \max_{\cJ \subset [N], |\cJ| = n} \quad \sum_{i \in \cI} \sum_{j \in \cJ} X_{ij}.
\eeq
One rejects $H_0$ when either of the two values go beyond some corresponding pre-defined thresholds. Note that in order to make the scan statistic \eqref{scan} work, exact knowledge of the size of $\cI$ and $\cJ$ is required, which here is $m$ and $n$. When $(m,n)$ are unknown, a Bonferroni testing procedure is used - that is, test all combinations of $(m,n)$ of interest, and reject $H_0$ when one of the tests indicates to reject. 

The method of \citep{butucea2013detection} relies on parametric assumptions, therefore \citep{arias2017distribution} considered calibrate the p-value by permutation. Since the permutation test is based on \eqref{scan}, full knowledge of $(m,n)$ is necessary. That is, it can only distinguish between $H_0$ and 
\begin{center}
$H_1(m,n)$: There exists $\cI\subset [M]$ and $\cJ \subset [N]$ such that $X_{ij}$ are larger than usual when $(i,j) \in \cI \times \cJ$, where $|\cI| = m$ and $|\cJ| =n$.
\end{center} 

In this paper we adapt the permutation test framework, and develop a Bonferroni testing procedure in order to deal with the case when $(m,n)$ are unknown. 

\begin{con}
We develop a Bonferroni testing procedure based on the permutation test by \citep{arias2017distribution}. We show that the testing procedure (at the first order) is asymptotically as powerful as the test by \citep{arias2017distribution} and \citep{butucea2013detection}. This is analyzed and proved under some standard exponential family parameter assumptions.
\end{con}

As permutation test is computationally hard to calibrate (usually the total number of permutations are increasing exponentially with the sample size), in practice a Monte Carlo calibration method is executed. This requires independent sampling from the group of permutation patterns for a decent number of times. However still, due to the computational complexity of \eqref{scan}, performing permutation tests on all combinations of $(m,n)$ could be extremely consuming in time and computational power. We construct a subset of $[M] \times [N]$, such that by performing permutation test on all $(m,n)$ inside this subset, we can still detect the existence of the elevated submatrix with no sacrifice of the power on the first order. 

\begin{con}
We propose a power-preserving fast test based on the Bonferroni permutation test, with the Bonferroni procedures working on a proper approximate net of $[M] \times [N]$. We show that this test is as powerful as the Bonferroni permutation test on all pairs of $(m,n)$ in $[M] \times [N]$. This is also analyzed and proved under some standard exponential family parameter assumptions.
\end{con}

\subsection{More Related work}

There are works focusing on the localization theory side of this problem, studying the existence of consistent estimators of the elevated matrix under some parametric setup \citep{kolar2011minimax,chen2016statistical, hajek2015information, hajek2015submatrix}. With the realization of the importance of computational efficiency, there is a stream of work line considering the trade-off between statistical and computational power \citep{balakrishnan2011statistical,cai2015computational, chen2016statistical, ma2015computational}. Several computationally efficient (with polynomial computational complexity with respect to data size) submatrix localization algorithm are proposed, including convex optimization \citep{chen2016statistical, chi2016convex}  and  spectral method \citep{cai2015computational}.

Another active research line worth mentioning here is on the stochastic block model (SBM). In the setup of SBM the observation is a graph, with edges independently connected. See \citep{holland1983stochastic} for a detailed introduction. The detection problem is to distinguish between an Erd\H{o}s-R\'{e}nyi graph and a graph with groups of nodes which nodes are more likely to connect within groups compared to across groups. The localization problem is to cluster the nodes by the closeness of their connection.

If the adjacency matrix of the graph is considered, the problem shares many properties with submatrix detection and localization (the adjacency matrix is symmetric with independent upper triangle nodes, compared with submatrix localization problem). There are works considering the existence of consistency detectors \citep{zhang2016minimax}, the existence of consistency clustering methods \citep{mossel2015consistency}, semi-definite programming \citep{chen2016statistical, abbe2016exact}, and spectral methods \citep{chaudhuri2012spectral, mcsherry2001spectral}.

\subsection{Content}
The paper is arranged as follows. \secref{bound} introduces the parametric setup, the detection boundaries set up by \citep{butucea2013detection}, as well as the permutation test by \citep{arias2017distribution}. \secref{bonf} describes the Bonferroni-type testing procedure as well as its theoretical property. \secref{algo} shows the construction of an approximation net in order to speed up the testing process, and the associated theoretical results. The numerical experiments are in \secref{numerics}.
\section{The detection boundaries}
\label{sec:bound}

We establish the minimax framework corresponding to the hypothesis testing problem raised in \secref{intro}. First we introduce the one-parameter natural exponential family considered throughout this paper, which is also the parametric family used in the analysis of \citep{arias2015distribution, butucea2013detection, arias2017distribution}. 
To define such a distribution family, first consider a distribution $\nu$ with mean zero and variance $1$, and assume that its moment generating function $\phi(\theta) < \infty$ for some $\theta>0$. Let $\theta_* = \sup\{ \theta: \phi(\theta) < \infty\}$ , and the family is parameterized by $\theta \in [0, \theta_*)$ with density function 
\beq \label{exp}
f_\theta(x) = \exp\{x\theta - \log \phi (\theta) \}.
\eeq 
The density function is with respect to $\nu$. Note that when $\theta=0$, $f_0(x) = 1$ which corresponds to $\nu$. By varying the choice of $\nu$, this parametric model includes several common parametric models such as normal family ($\nu = \mathcal{N}(0,1)$), Poisson family ($\nu = \text{Pois}(1) - 1 $) and Rademacher family ($\nu = 2 \text{Rade}(0.5) -1$).

An important property of $f_\theta$ is the stochastic monotonicity with respect to $\theta$. This enables us to model the previous hypothesis testing problem with $\nu$ acting as noise distribution and $f_\theta$ as the distribution of those unusually large entries. In details, we consider 
\begin{center}
$H_0$: $X_{ij}$ are IID following $\nu$ for all $i\in [M]$ and $j\in [N]$ \\ 
vs\\
$H_1$: There exists $\cI\subset [M]$ and $\cJ \subset [N]$ such that $X_{ij} \sim f_{\theta_{ij}}, \theta_{ij} > \theta_\ddag > 0
$ when $(i,j) \in \cI \times \cJ$.
\end{center}
Here the parameter $\theta_\ddag$ is the lower bound for all the $\theta_{ij}$ inside the raised submatrix, and is acting as the role of signal-noise ratio. 
In this context, if the potential submatrix's size $(m,n)$ is known, the corresponding $H_1(m,n)$ is described as follows.

\begin{center}
$H_1(m,n)$: There exists $\cI\subset [M]$ and $\cJ \subset [N]$ such that $X_{ij} \sim f_{\theta_{ij}}, \theta_{ij} > \theta_\ddag > 0
$ when $(i,j) \in \cI \times \cJ$, where $|\cI| = m$ and $|\cJ| =n$.
\end{center}

In ordert to perform the hypothesis testing task between $H_0$ and $H_1(m,n)$ under this parameterization framework, \cite{butucea2013detection} developed the following.

\begin{thm}[\cite{butucea2013detection}] \label{thm:bi}
Consider an exponential model as described in \eqref{exp}, with $\nu$
having finite fourth moment. Assume that 
\beq  \label{cond}
M,N,m,n\to \infty, \quad m = o(M), n= o(N), \quad \frac{\log(M\vee N)}{m \wedge n} \to 0.
\eeq
For any $\alpha > 0$ fixed, the sum test based on \eqref{sum} is asymptotically powerful if
\beq 
\frac{\theta_\ddag mn}{\sqrt{MN} } \to \infty,
\eeq
and the scan test based on \eqref{scan}, is asymptotically powerful if 
\beq \label{condmain}
\lim \inf \frac{\theta_\ddag \sqrt{mn}}{\sqrt{2(m\log(M/m) + n\log(N/n))}} > 1.
\eeq
Conversely, if $m = O(n)$ and $\log M = O(\log N)$, and 
\beq \label{lowerbound}
\frac{\theta_\ddag mn}{\sqrt{MN} } \to 0, \quad \lim \inf \frac{\theta_\ddag \sqrt{mn}}{\sqrt{2(m\log(M/m) + n\log(N/n))}} < 1,
\eeq 
any test with level $\alpha$ has limiting power at most $\alpha$. 
\end{thm}

The test rejects $H_0$ when the sum test statistic \eqref{sum} or scan test statistic \eqref{scan} is larger than some pre-defined threshold. While these tests heavily depend on parameterization, distribution-free methods such as permutation test are proposed to tackle this problem. A permutation test based on the scan statistic \eqref{scan} is analyzed by \citep{arias2017distribution}. Here a permutation pattern on set $[M] \times [N]$ is denoted by $\pi$, and define
\begin{itemize}
\item $\Pi_1$ as the collection of permutations that entries are permuted within their row;
\item $\Pi_2$ as the collection of all permutations.
\end{itemize}
To illustrate the difference of $\Pi_1$ and $\Pi_2$, let $\bX_\pi = \{X_{\pi(i,j)}\}$, consider 
\beqn 
\bX = \begin{pmatrix} 
1 & 2 &3 &4\\
5 & 6 &7 &8 \\
9 &10 &11 &12
\end{pmatrix}
\eeqn
and we will have 
\beqn 
\bX_{\pi_1} = \begin{pmatrix} 
1 & 4 &2 &3\\
6 & 8 &7 &5 \\
10 &12 &9 &11
\end{pmatrix},\bX_{\pi_2} = \begin{pmatrix} 
10 & 8 &2 &9\\
6 & 3 &11 &1 \\
12 &4 &5 &7
\end{pmatrix} 
\eeqn
for some
$\pi_1 \in \Pi_1, \pi_2 \in \Pi_2$. 
We calculate the scan statistic on the permuted data $\bX_\pi$, and compare to the one from the original data $\bX$. The permutation $p$-value is defined as 
\beq \label{perm}
\mathfrak{P}_{m,n}(\bX) = \frac{|\{\pi \in \Pi : \scan_{m,n}(\bX_\pi ) \geq \scan_{m,n}(\bX) \}| }{|\Pi|}.
\eeq
Here $\Pi =\Pi_1$ or $\Pi_2$. The permutation test rejects the null hypothesis when $\mathfrak{P}_{m,n}(\bX)$ is smaller than the pre-determined level. This test has the following property.
\begin{thm}[\cite{arias2017distribution}, Theorem 2] \label{thm:acl}
Consider an exponential model as described in \eqref{exp}, assume \eqref{cond} and 
\beq  \label{extracond}
\log^3 (M \vee N) / (m \wedge n)  \to  0
\eeq
and that either (i) $\nu$ has support bounded from above, or (ii) $\max_{i,j} \theta_{ij} \leq \bar{\theta}$ for some $\bar{\theta} < \theta_*$ fixed.
Additionally if $\Pi_i = \Pi_1$, we require that $\phi(\theta) < \infty$ for some
$\theta < 0$. Then the permutation scan test based on \eqref{perm}, at any fixed level $\alpha > 0$, has limiting power $1$ 
when \eqref{condmain} holds.
\end{thm}

This theorem shows that permutation test has the same first-order asymptotic power compared to the parametric test described in \thmref{bi}, under some extra mild conditions. 

\section{Bonferroni permutation test}
\label{sec:bonf}

We now consider the case where the size of the submatrix to be detected is not specified. Recall the definition of $H_1$ and $H_1(m,n)$ in \secref{bound}. By realizing the fact that $H_1$ is true if and only if there exists some $(m,n)$ such that $H_1(m,n)$ is true, we can perform test on $H_1(m,n)$ for all pairs of $(m,n)$, and use Bonferroni correction in order to control the type I error. 

We adapt the distribution-free permutation test of \cite{arias2017distribution} here. For each pair of $(m,n)$, calculate $\mathfrak{P}_{m,n}(\bX)$, and calculate the final Bonferroni corrected $p$-value as 
\beq \label{bonf}
\mathfrak{P}(\bX) = \min (MN \min_{m,n} \mathfrak{P}_{m,n} (\bX),1).
\eeq
One rejects $H_0$ when $\mathfrak{P}(\bX)$ is less than some pre-determined level $\alpha$. Due to the property of Bonferroni type of tests, this test has level $\alpha$ if each test concerning $H_0$ and $H_1(m,n)$ has level $\alpha / MN$, which is a proved fact in \citep{arias2017distribution}, regardless of the dependencies between tests. 

Being a conservative method in multiple testing, Bonferroni method usually loses statistical power in exchange for controlling the family wise error rate. However in some cases (for example, \citep{arias2014community, butucea2013detection}), the Bonferroni procedure achieves the same first order asymptotic power as the scan test without knowledge of the submatrix size. We illustrate the same phenomenon in this case. 

\begin{thm} \label{thm:bonf}
Consider an exponential model as described in \eqref{exp}. Assume that there exists a pair of $(m,n)$ such that $H_1(m,n)$ is true, and all the assumptions in \thmref{acl} are satisfied, which are:
\beq \label{condbonf}
M,N,m,n\to \infty, \quad m = o(M), n= o(N), \quad \frac{\log(M\vee N)}{m \wedge n} \to 0, \quad \log^3 (M \vee N) / (m \wedge n)  \to  0
\eeq
and that either (i) $\nu$ has support bounded from above, or (ii) $\max_{i,j} \theta_{ij} \leq \bar{\theta}$ for some $\bar{\theta} < \theta_*$ fixed.
Additionally if $\Pi_i = \Pi_1$, we require that $\phi(\theta) < \infty$ for some
$\theta < 0$. 
Then 
\beq 
\mathfrak{P}(\bX) \to 0
\eeq
in probability.
\end{thm}

The interpretation of this theorem is as follows. Assume we are under $H_1(m,n)$. With full knowledge of $(m,n)$, if the permutation test by \eqref{perm} can successfully detect the submatrix with high probability, the test by \eqref{bonf} will reject $H_0$ with high probability, without any knowledge of $(m,n)$.  

Note that if \eqref{lowerbound} and its associated assumptions hold true, any test trying to distinguish $H_0$ and $H_1$ with level $\alpha$ will have limiting power at most $\alpha$, due to the fact that any test distinguishing $H_0$ versus $H_1$ can be used to distinguish $H_0$ and $H_1(m,n)$, thus having limiting power at most $\alpha$ thanks to \thmref{bi}. 

\section{A power-preserving fast test}
\label{sec:algo}

We build our testing framework on the permutation test, which is by nature a computationally intensive method. Calculation of scan statistic \eqref{scan} is NP-hard, and the total number of permutations $|\Pi|$ will skyrocket when the number of data increases. In practice the scan statistic is calculated by LAS algorithm proposed by \cite{shabalin2009finding}, as did in previous literature \citep{butucea2013detection, arias2017distribution}, and the permutation test is done by Monte-Carlo sampling. In detail, a large number $B$ is fixed, and permutations $\{ \pi_1, \ldots, \pi_B\}$ is the IID sample from uniform distribution on $|\Pi_i|$. Then $\mathfrak{P}_{m,n}$ is approximated by 
\beq \label{mc}
\hat{\mathfrak{P}}_{m,n} (\bX)= \frac{|\{i:  \scan_{m,n}(\bX_{\pi_i} ) \geq \scan_{m,n}(\bX) \}| +1 }{B + 1}.
\eeq

As illustrated in \citep{butucea2013detection, arias2017distribution}, the calculation of the scan statistic \eqref{scan} is already difficult,
even if the size $(m,n)$ is known. Now with the submatrix size unknown, we have added
the difficulty since the scan statistic under all possible combinations of $(m,n)$, will be
calculated during the Bonferroni process. In principle, the Bonferroni method requires
going over all submatrix sizes, but we only scan a carefully chosen subset to lighten up
the computational burden. Inspired by \cite{arias2005near, arias2015distribution}, in which the authors scan for anomalous data interval within all possible intervals using a dyadic representation, we illustrate the construction of such subset on $[M] \times [N]$, and show that the first-order statistical power is preserved.

\subsection{An approximation net}

The subset of $[M] \times [N]$ we are going to construct in order to approximate the elements in $[M] \times [N]$ is called an approximate net. We first construct one-dimensional approximate net on $[M]$. 

We start by the following definition.
\begin{Def}\label{def:binary}
A binary expansion of an integer $c$ is a sequence $\{a_i(c)\}$ with $ a_i(c) \in \{0,1\}$ and $i \leq \lfloor \log_2 c \rfloor$, such that
\beq
c = \sum_{i = 0}^{\lfloor \log_2 c \rfloor} a_i(c) 2^i.
\eeq
\end{Def}
After representing an integer in the binary numeral system, one may approximate this integer by keeping the first $k$ digits of its binary expansion.
\begin{Def}\label{def:kbinary}
	A $k$-binary approximation of an integer $c$ is the interger $c'$ such that
\beq
	a_i(c') = a_i(c) \1 \{i \geq \lfloor \log_2 c \rfloor - k+1\}.
\eeq
\end{Def}
To find the $k$-binary approximation of $c$, represent $c$ in its binary expansion, keep the first $k$ digits, and shrink the rest to zero. Finally calculate $c'$ by the formula in \defref{binary}. It captures the main part of the integer
and the difference could be controlled by $k$. The method is closed related to the binary tree representation of integers, where an integer is represented to the root at distance $k$ from the first non-zero node. The following lemma gives an upper bound of the difference rate.

\begin{lem} \label{lem:binary}
	If $c'$ is a $k$-binary approximation of $c$, then
\beq
\frac{c - c'}{c} \leq 2^{1-k}.
\eeq
\end{lem}
Note that the difference rate is only associated with $k$, but not with the value of $c$. Therefore if we apply $k$-binary approximation to a collection of integers, the difference rate will be controlled \emph{uniformly} among all the integers in the collection by the choice of $k$. 

Now we construct the approximation net of $[M]$ based on $k$-binary approximation.
\begin{Def} \label{def:appnet}
	An approximation net $S_k(M)$, based on $k$-binary approximation, of set $[M]$, is defined as
\beq
	S_k(M)=\{c' : \text{$c'$ is a $k$-binary approximation of some $c \in [M]$} \}
\eeq
\end{Def}
To better illustrate the approximation net, \tabref{approx_net_1024} shows the approximation net of $1024$ with $k=3$.

\begin{table}
\begin{center}
\caption{Approximation net for $1024$}
\label{tab:approx_net_1024}
\vspace{2pt}
\begin{tabular}{clclclclclclc}
\hline
Binary & Decimal & Binary & Decimal & Binary & Decimal\\
\hline
10000000000 & 1024 & 00010000000 & 128 & 00000010000 & 16   \\
01110000000 & 896 & 00001110000 & 112 & 00000001110 & 14   \\
01100000000 & 768 & 00001100000 & 96  & 00000001100 & 12   \\
01010000000 & 640 & 00001010000 & 80  & 00000001010 & 10   \\
01000000000 & 512 & 00001000000 & 64  & 00000001000 & 8    \\
00111000000 & 448 & 00000111000 & 56 &  00000000111 & 7    \\
00110000000 & 384 & 00000110000 & 48 &  00000000110 & 6    \\
00101000000 & 320 & 00000101000 & 40 &  00000000101 & 5    \\
00100000000 & 256 & 00000100000 & 32 &  00000000100 & 4    \\
00011100000 & 224 & 00000011100 & 28 &  00000000011 & 3    \\
00011000000 & 192 & 00000011000 & 24 &  00000000010 & 2    \\
00010100000 & 160 & 00000010100 & 20 &  00000000001 & 1   \\

\hline
\end{tabular}
\end{center}
\end{table}

The cardinality of $S_k(M)$ can be much less than $M$ if $k$ is chosen properly. For example, set $k = \log_2 \log_2 M + 1$ and it can be shown that 
$|S_k(M)| = O((\log_2 M)^2)$. Note that in this case $k \to \infty$ when $M \to \infty$, and by \lemref{binary} we know that for every $c \in [M]$, there exists some $c' \in S_k(M)$ such that $c' = (1 + o(1)) c$, and $o(1)$ is uniform among $[M]$.  

\subsection{Test power under approximation nets}

Based on the one-dimensional approximation net defined in \defref{appnet}, we can similarly extend the idea to sets of two-dimensional integer pairs. We perform the Bonferroni-type testing procedure on $S_{k_M}(M) \times S_{k_N}(N)$, instead of $[M] \times [N]$. In detail, we use the following Bonferroni corrected $p$-value:
\beq \label{accbonf}
\mathfrak{P}_{k_M, k_N}(\bX) = \min ( |S_{k_M}(M)| |S_{k_N}(N)| \min_{(s,t) \in S_{k_M}(M) \times S_{k_N}(N)} \mathfrak{P}_{s,t}(\bX), 1 ).
\eeq
The idea is to use the property of approximation net to eliminate a significant portion of calculation by reducing the scanning region of the Bonferroni process, while keeping the accuracy through choosing a proper pair of $(k_M, k_N)$. Assume we are under $H_1(m,n)$. When setting $k_M, k_N \to \infty$, there is a pair of $(m', n')\in S_{k_M}(M) \times S_{k_N}(N)$ that is close enough to $(m,n)$, and $\mathfrak{P}_{m',n'}(\bX)$ will converge to zero fast enough such that brings the Bonferroni corrected $p$-value to zero as well. 

The following theorem describes the asymptotic power of the Bonferroni test on the approximate net.
\begin{thm} \label{thm:bonf2}
Assume all the assumptions in \thmref{bonf} hold (specifically, \eqref{condbonf} and (i) and (ii) following it). Further set $k_M, k_N \to \infty$. Then
\beq
\mathfrak{P}_{k_M, k_N}(\bX) \to 0
\eeq
in probability.
\end{thm}

\section{Numerical experiments}
\label{sec:numerics}



We use simulations to verify our theoretical findings.\footnote{The code used in this section is available in \url{https://github.com/nozoeli/bonferroniSubmatrix}} The main purpose of the simulation is to illustrate the proposed test's behavior under the alternative hypothesis. We adapt the simulation setup in \cite{arias2017distribution} due to the similarity of the research problem and simulation purpose. 

The data $\bX$ is generated with IID samples from $\nu$ for all entries except $[m] \times [n]$, where the entries are IID from $f_\theta$ for some $\theta > 0$. Then the data goes through the proposed test progress and the $p$-value is recorded. The above process is repeated $100$ times for fixed $\theta$ in order to see the stochastic behavior of the $p$-values. 

To see how the theory works, we define $\theta_{\text{crit}}$ as follows.
\beq
\theta_{\text{crit}} = \sqrt{\frac{2(m\log (M/m) + n \log (N/n))}{mn}}.
\eeq
And we slowly increase $\theta$ from $0.625 \times \theta_{\text{crit}}$ to $1.5 \times \theta_{\text{crit}}$, with step size $0.125\times\theta_{\text{crit}}$. This is aim to examine the behavior of the test around $\theta = \theta_{\text{crit}}$, which is claimed to be the convergence threshold from our theory. 

The permutation test is calibrated by Monte-Carlo which is described in \eqref{mc}, with $B = 500$. As mentioned in \secref{algo}, the calculation of $\scan_{m,n}(\bX)$ is NP-hard in theory. LAS algorithm from \cite{shabalin2009finding} is an approximating algorithm to calculate the scan statistic, however due to its hill-climbing optimization process, it suffers from being stuck inside local minimums. We re-initiate the LAS algorithm several times with random initialization, and return the result with the largest output, in order to prevent being stuck at local minimums. Both permutation methods are examined in the analysis. The data size is set as $(M,N) = (200,100)$ and we examine $2$ anomaly sizes, namely $(m,n) = (10,15)$ and $(30,10)$. For the approximate net, we set $(k_M, k_N) = (\lfloor \log_2(\log_2(M)) \rfloor , \lfloor \log_2(\log_2(N)) \rfloor)$

\begin{figure}%
	\centering
	\subfloat{{\includegraphics[width=0.8\textwidth]{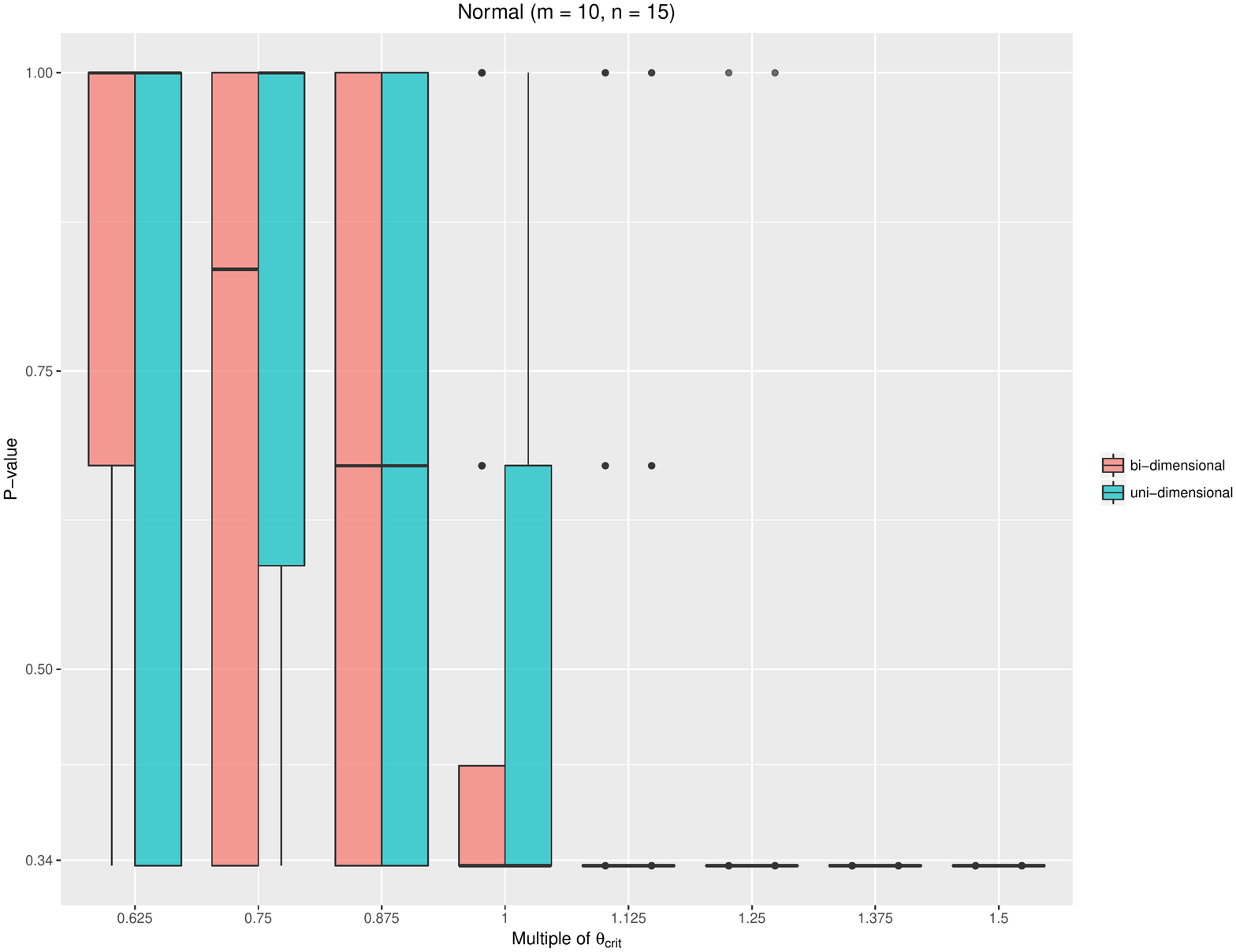} }}\\%
	\subfloat{{\includegraphics[width=0.8\textwidth]{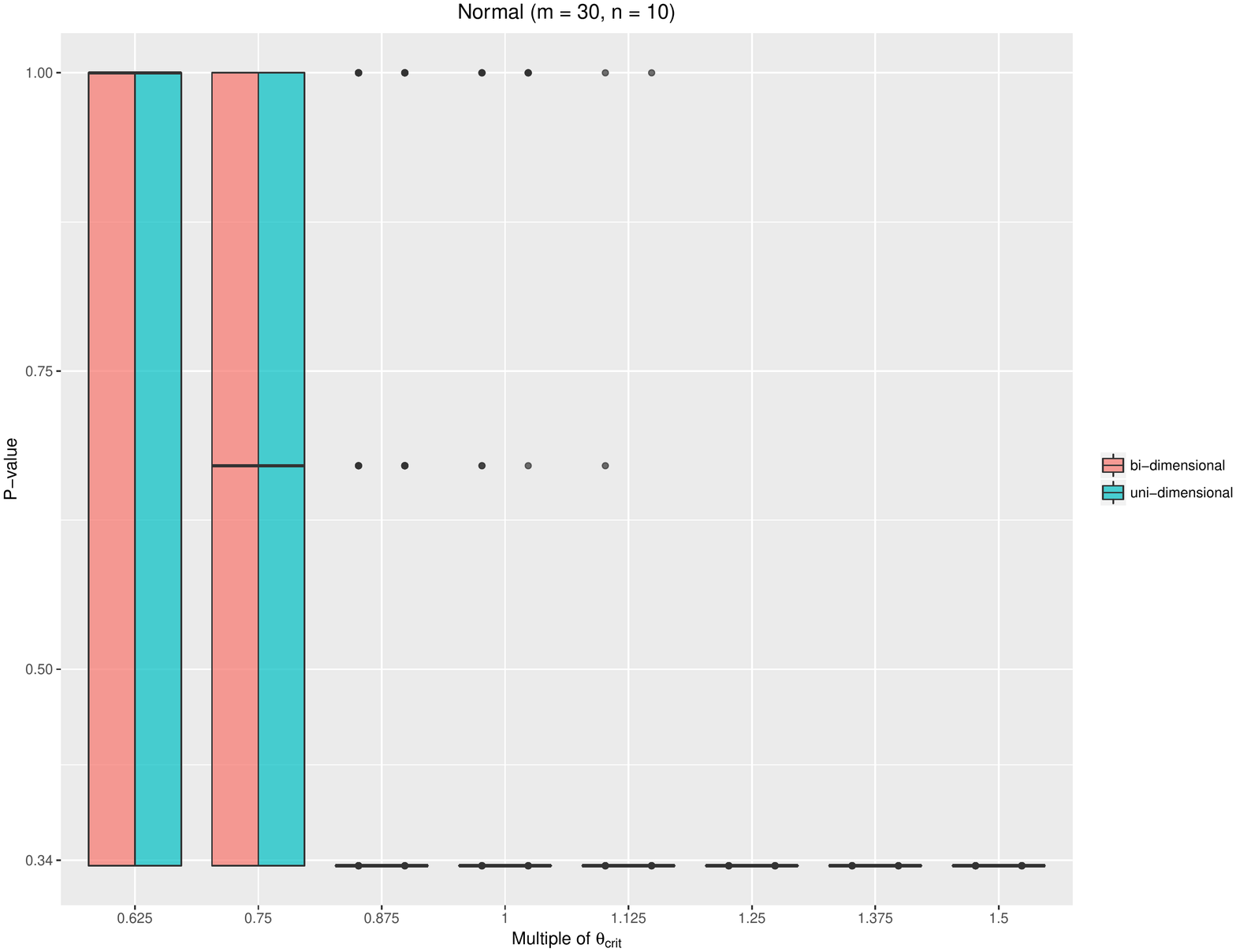} }}%
	\caption{$p$-values of proposed test in the normal model}%
	\label{fig:norm}%
\end{figure}

\begin{figure}%
	\centering
	\subfloat{{\includegraphics[width=0.8\textwidth]{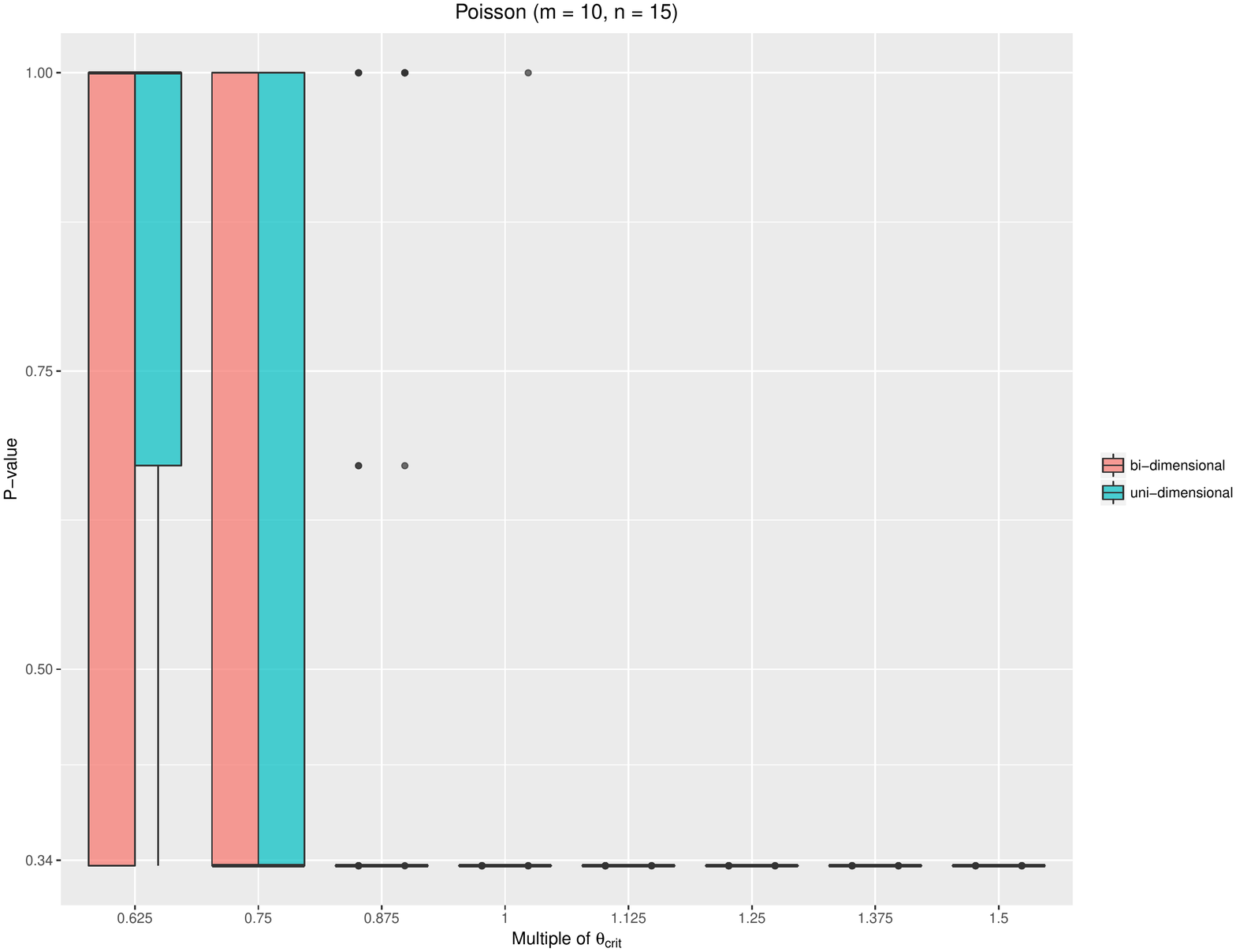} }}\\%
	\subfloat{{\includegraphics[width=0.8\textwidth]{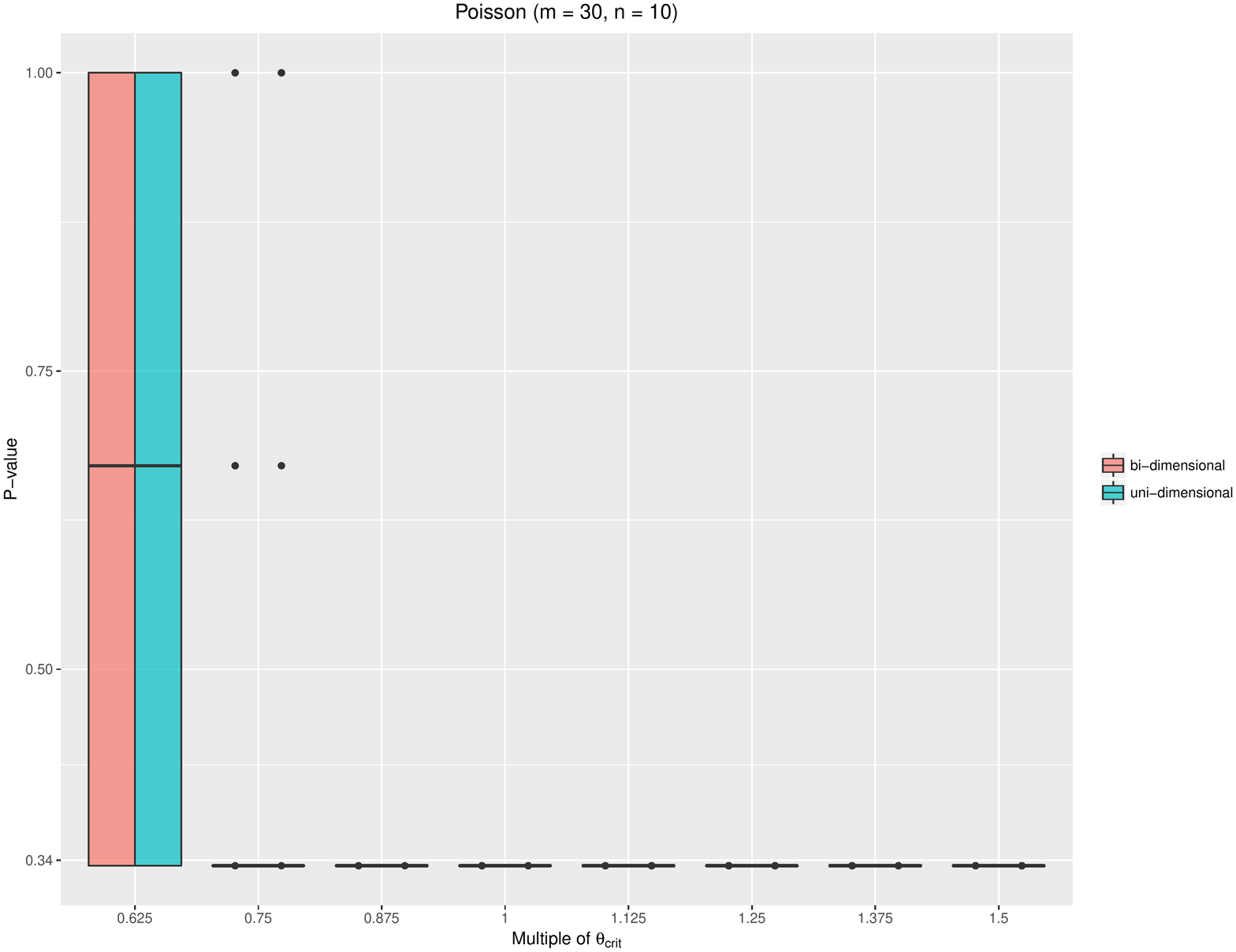} }}%
	\caption{$p$-values of proposed test in the Poisson model}%
	\label{fig:poi}%
\end{figure}

We choose two representative distributions as $\nu$, which is standard normal and centralized $Pois(1)$. The corresponding $f_\theta$ is $\mathcal{N} (\theta, 1)$ and $Pois(e^\theta) - 1$. From Figure \ref{fig:norm} and \ref{fig:poi} we see that when $\theta \leq \theta_{\text{crit}}$, the $p$-values are generally close to $1$ due to the conservative property of the Bonferroni methods. However, the $p$-value shrinks to simulation lower bound at 
\beq
|S_k(M)| \times |S_k(N)|  / (B + 1) = 14 \times 12 / (500 + 1) \approx 0.34\eeq 
quickly once $\theta$ passes $\theta_{\text{crit}}$, in both permutation patterns and distribution setups. This is because of the fact that all the permutation $p$-values are larger than or equal to $1/(B+1)$, which makes the Bonferroni corrections no less than $|S_k(M)| \times |S_k(N)|  / (B + 1)$. Although reaching the lower bound indicates that no permutation during the test generates the test statistic larger than the one from the original data, it is not very persuasive in the sense to show that the $p$-values are converging to $0$. 
In order to verify this, we need a larger number of Monte Carlo simulations, namely a larger $B$, such that $|S_k(M)| \times |S_k(N)|  / (B + 1)$ is close to $0$.

\begin{figure}%
	\centering
	\subfloat{{\includegraphics[width=0.8\textwidth]{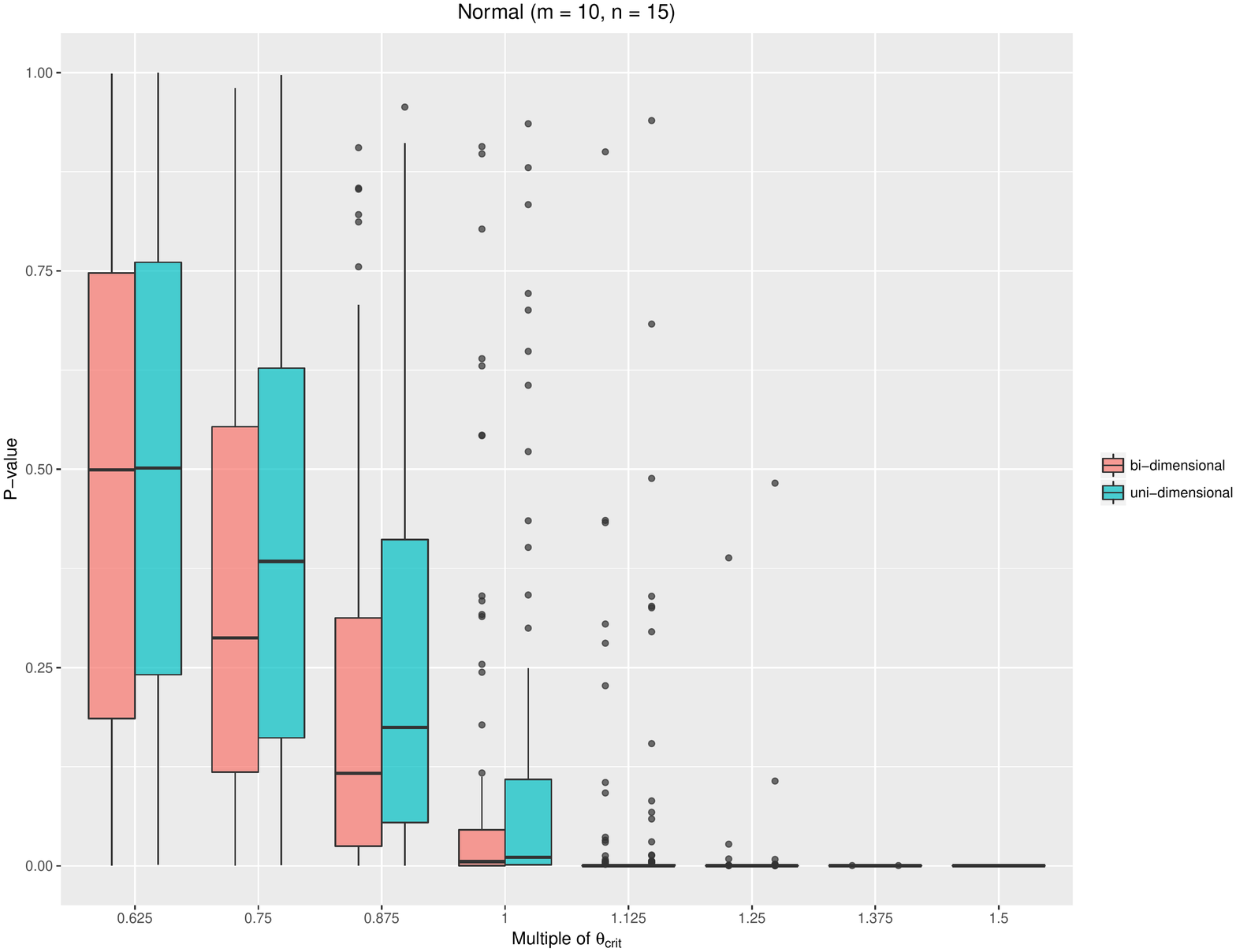} }}\\%
	\subfloat{{\includegraphics[width=0.8\textwidth]{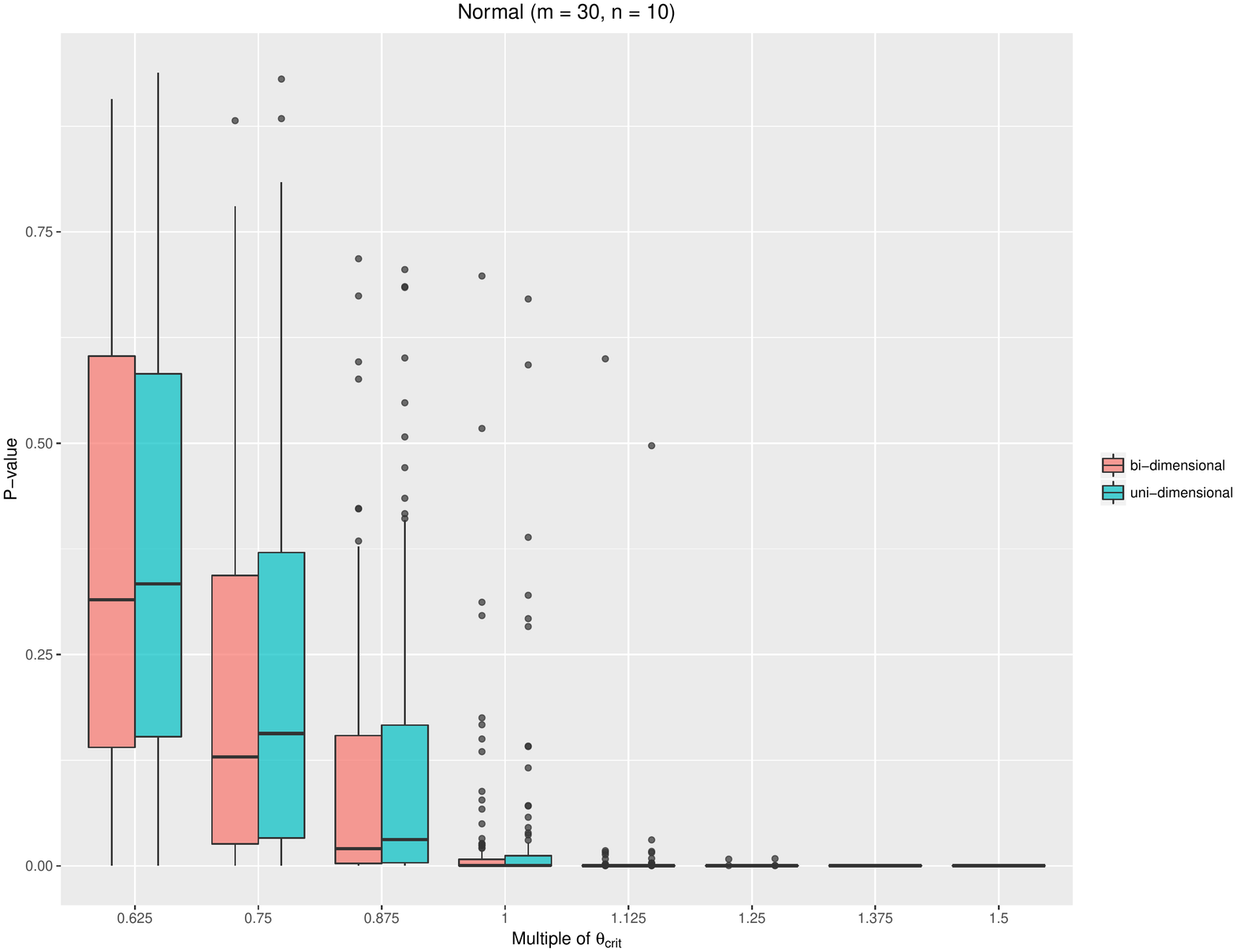} }}%
	\caption{Upper bounds of $p$-values of proposed test in the normal model}%
	\label{fig:norm_upper}%
\end{figure}

\begin{figure}%
	\centering
	\subfloat{{\includegraphics[width=0.8\textwidth]{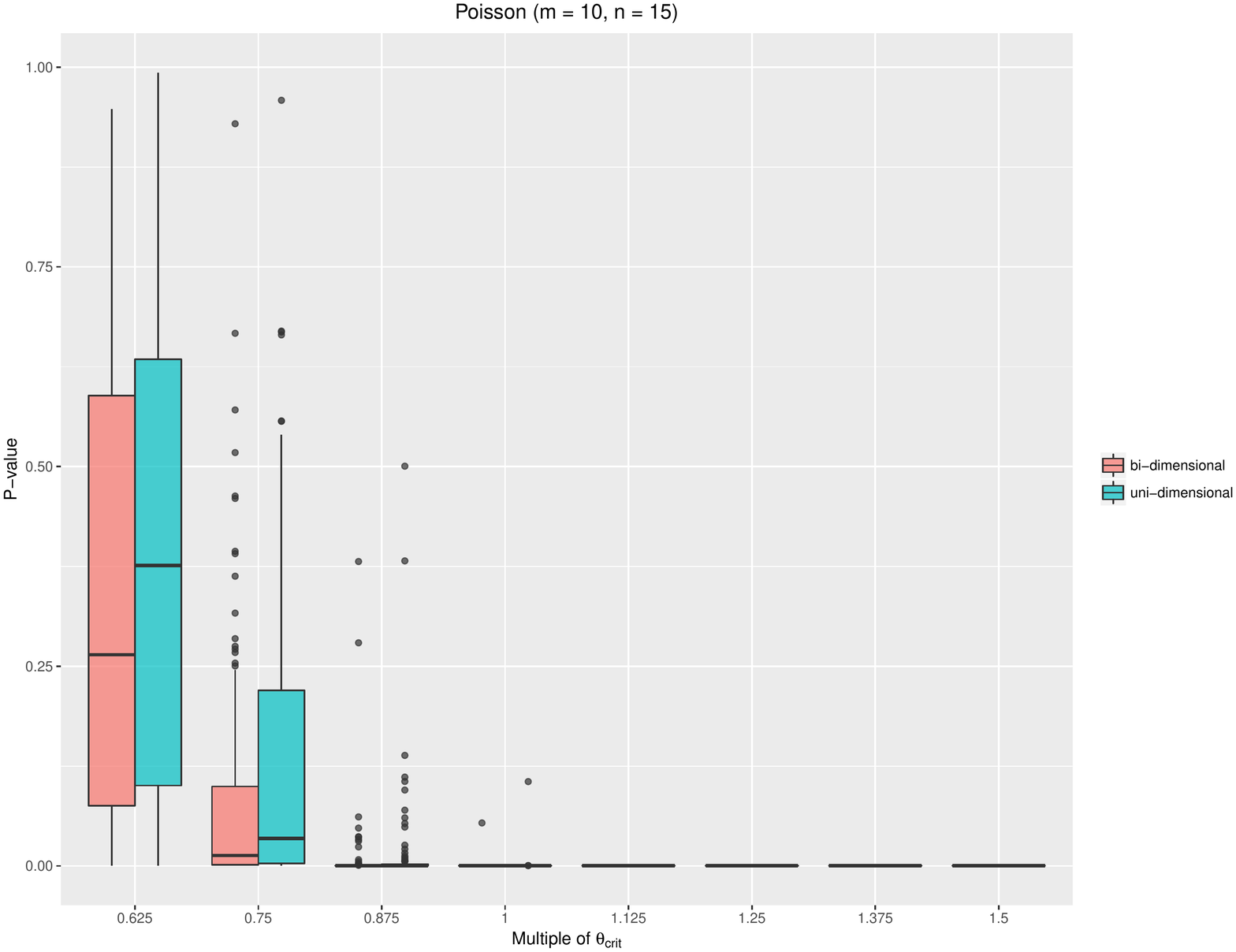} }}\\%
	\subfloat{{\includegraphics[width=0.8\textwidth]{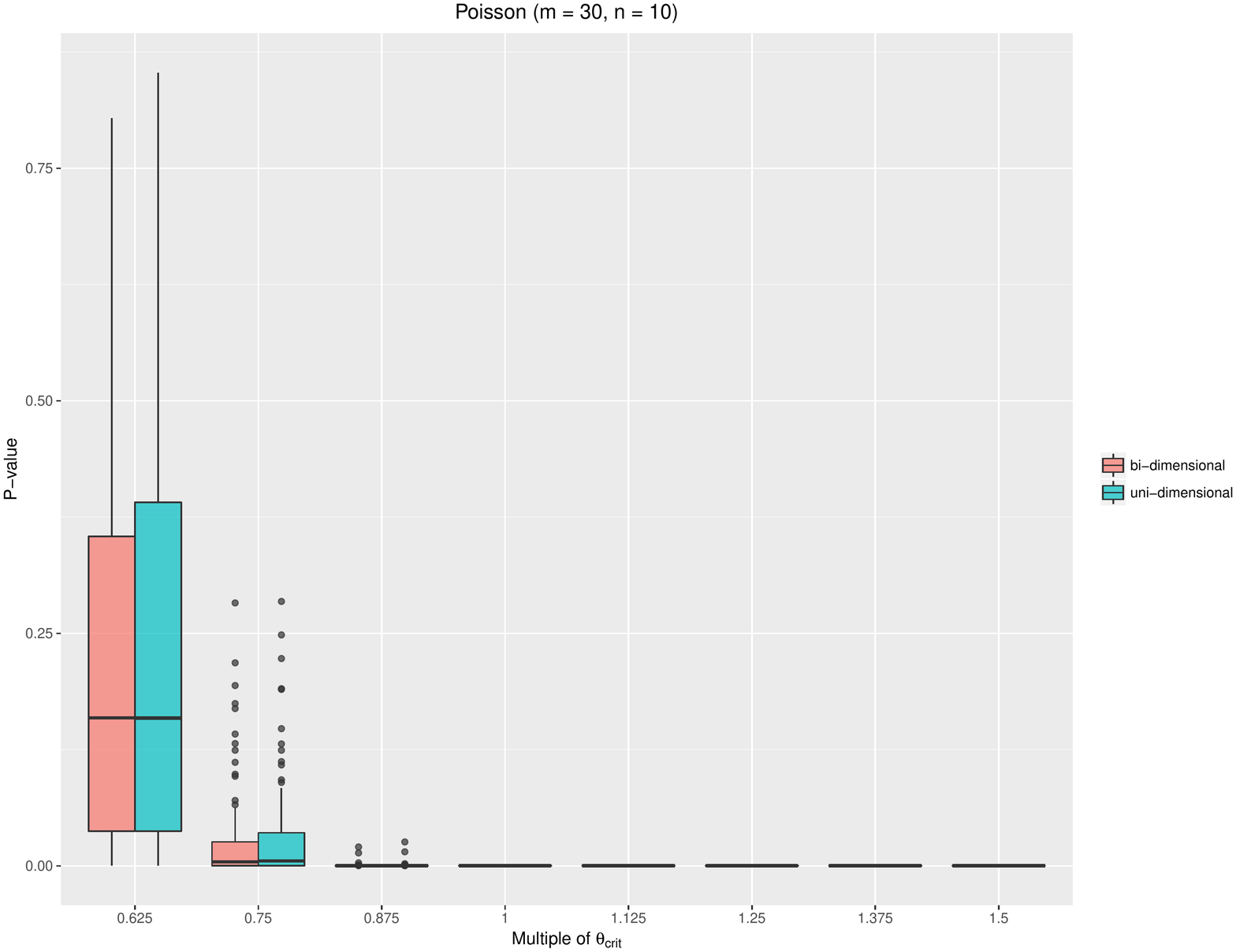} }}%
	\caption{Upper bounds of $p$-values of proposed test in the Poisson model}%
	\label{fig:poi_upper}%
\end{figure}

We perform an additional experiment by setting $B = 5000$, and focus on the pair of $(m', n')\in S_{k_M}(M) \times S_{k_N}(N)$, such that $m' = \min \{ p \in S_{k_M}(M) | p > m \}$ and $n' = \min \{ q \in S_{k_N}(N) | q > n \}$. The provides us an upper bound of the Bonferroni permutation $p$-value, since by \eqref{accbonf},
\beq 
\mathfrak{P}_{k_M, k_N}(\bX) \leq \min ( |S_{k_M}(M)| |S_{k_N}(N)|  \mathfrak{P}_{m',n'}(\bX), 1 ).
\eeq
We illustrate the relationship between this upper bound of the Bonferroni permutation $p$-value and the signal level. \figref{norm_upper} and \figref{poi_upper} shows that the upper bound converges to $0$ once the signal level $\theta$ goes beyond $\theta_{\text{crit}}$, and this confirms our theoretical findings. 

\section{Technical proofs}
\label{sec:proof}

In this section we present the proofs to the theorems in the previous sections. Unless separately declared, the ineqalities holds with high probability (or with probability $1 -o(1)$ as data matrix size $M,N \to \infty$). 

The following lemma is at the center of our argument. 

\begin{lem}[\citep{arias2015distribution} Lemma 2, Bernstein's inequality for sampling without replacement]
	\label{lem:bernstein}
	Let $(Z_1, \dots, Z_m)$ be obtained by sampling without replacement from a given a set of real numbers $\{z_1, \dots, z_J\} \subset \bbR$.  Define 
	$z_{\rm max} = \max_j z_j$, $\bar z = \frac1J \sum_j z_j$, and $\sigma_z^2 = \frac1J \sum_j (z_j - \bar z)^2$.
	Then the sample mean $\bar Z = \frac1m \sum_i Z_i$ satisfies 
	\[
	\P\big(\bar Z \ge \bar z + t\big) 
	\le \exp\Bigg[- \frac{m t^2}{2 \sigma_z^2 + \frac23 (z_{\rm max} - \bar z) t}\Bigg], \quad \forall t \ge 0.
	\]
\end{lem}

The lemma is a result from \cite{serfling1974probability}. We also refer the reader to \citep{bardenet2015concentration, boucheron2013concentration} for further details of the lemma.

\subsection{Proof of \thmref{bonf}}  \label{sec:proofbonf}
	We first show that the test has level $\alpha$. To show this, we fix a pair of $(m,n)$ and show that $\P (\mathfrak{P}_{m,n}(X) \leq \alpha)\leq \alpha$. By the standard argument on the level of Bonferroni test, we will finish the proof on the level. Assuming the null is true, $\textbf{scan}_{m,n}(\pi(X))$ has the same distribution with $\textbf{scan}_{m,n}(X)$ under either permutation methods. Therefore we define 
	$$
	T_k = \textbf{scan}_{m,n}(\pi_k (X)), k \in [(MN)!],
	$$ 
	and assume $T_{k_0} = \textbf{scan}_{m,n}(X)$, then $\text{rank} (T_{k_0})$ is uniformly distributed on $[(MN)!]$ (if the ties are broken randomly). We have 
	$$
	\P (\mathfrak{P}_{m,n}(X) \leq \alpha) \leq \P (\text{rank} (T_{i_0}) \leq \lfloor \alpha (MN!) \rfloor ) \leq \frac{\lfloor \alpha (MN!) \rfloor}{MN!} \leq \alpha.
	$$
	
	Therefore all we need to show is that the p-value tends to zero under the alternative. 
	
	\subsubsection{Bidimensional permutation}
	
	We start with considering bidimensional permutation. Starting from here we let $S^*$ be the true anomaly submatrix, and $(m,n)$ as its row and column size. For a submatrix index set $S$, let $Y_S(X) = \sum_{(i,j) \in S} X_{ij}$. Now fix some $S \in \mathbb{S}_{mn}$ and uniformly sample a permutation $\pi$. Here $\mathbb{S}_{mn}$ is the collection of row/column indexes of all submatrices of size $m\times n$. Also, if we conditional a realization of $\bX$, say $\bx$, by \lemref{bernstein}, we have
	$$
	\P \bigg(\frac{Y_{S} (\pi(\bx))}{mn} - \bar{\bx} > t\bigg) \leq \exp \Bigg( - \frac{t^2mn}{2\sigma^2_{\bx} + \frac{2}{3} (\bx_{\max} - \bar{\bx}) t} \Bigg),
	$$
	where $\bx_{\max}$ is the largest value in the data matrix $\bx$, and $\sigma_{\bx}^2$ is the sample variance of the data if $\bx$ is treated as a one dimensional data vector. Notice that the probability is on the permutation process. Apply a union bound, and set $t = \textbf{scan} (\bx)/mn - \bar{\bx}$, then we have
	\begin{eqnarray} 
	\mathfrak{P}(\bx)   \leq MN\mathfrak{P}_{m,n} (\bx)\leq MN|\mathbb{S}_{mn}| \exp\Bigg( - \frac{ mn(\textbf{scan}_{m,n} (\bx)/mn - \bar{\bx})^2}{2\sigma^2_{\bx} + \frac{2}{3} (\bx_{\max} - \bar{\bx}) ( \textbf{scan}_{m,n} (\bx)/mn - \bar{\bx})} \Bigg).
	\end{eqnarray}
	Note that this inequality holds for all the realizations of $\bX$, thus we allow $\bX$ to change and focusing on the quantity
	\beq \label{eq:bern}
	 MN|\mathbb{S}_{mn}| \exp\Bigg( - \frac{ mn(\textbf{scan}_{m,n} (\bX)/mn - \bar{\bX})^2}{2\sigma^2_{\bX} + \frac{2}{3} (\bX_{\max} - \bar{\bX}) ( \textbf{scan}_{m,n} (\bX)/mn - \bar{\bX})} \Bigg).
	 \eeq
	We show that this quantity is $o_P(1)$.
	
	We start with bounding $\bar{\bX}$ by rewriting $\bar{\bX}$ as
	\beq
	\bar{\bX} = \frac{\sum_{i,j} X_{ij}}{MN} = \frac{mn}{MN}\cdot  \frac{\sum_{(i,j)\in S^*} X_{ij}}{mn} + \frac{MN - mn}{MN} \cdot  \frac{\sum_{(i,j)\notin S^*}X_{ij} }{MN - mn}.
	\eeq
	With $\theta$ in the anomalous submatrix bounded from above, or the support of $f_\theta$ being bounded, the first term is $o_P (1)$ since $mn = o(MN)$. By Law of Large Numbers, the second term is $o_P(1)$ given that distribution $f_0$ has mean zero. So 
	\beq
	\bar{\bX} = o_P (1).
	\eeq
	
	Following the proof of Theorem 1 in \citep{arias2015distribution}, we can bound $\bX_{\max}$ with 
	\begin{eqnarray} \label{eq:xmax}
	\P	(\bX_{\max} - \bar{\bX} < \frac{3}{c} \log (MN)) \rightarrow 1.
	\end{eqnarray}
	with $c \in (0, \theta_* - \bar{\theta})$. Define the event $\mathcal{A} = \{\bX_{\max} - \bar{\bX} < \frac{3}{c} \log (MN)\}$. All the following arguments are conditional on $\mathcal{A}$. But since $\mathcal{A}$ happens with probability tending to $1$, all the conditional high probability events will happen unconditionally with high probability as well. 
	
	We do similar operations to bound $\sigma^2_\bX$ as follows.
	\begin{eqnarray}
		\sigma^2_{\bX} &=& \frac{1}{MN} \sum_{i,j} (X_{ij} - \bar{\bX})^2 \leq \frac{1}{MN} \sum_{i,j} X_{ij}^2 \\
		&=& \frac{mn}{MN}\cdot  \frac{\sum_{X_{ij}\in S^*} X_{ij}^2}{mn} + \frac{MN - mn}{MN} \cdot  \frac{\sum_{X_{ij}\notin S^*}X_{ij}^2 }{MN - mn}.
	\end{eqnarray}
	The first term is $o_P(1)$ if distribution $f_\theta$ has finite second moment, which can be derived from the assumption. The second term is $1 + o_P(1)$ by Law of Large Numbers and Slutsky's Lemma. Therefore
	\begin{eqnarray} \label{eq:var}
	\sigma_{\bX}^2 = 1 + o_P(1)
	\end{eqnarray}
	
	Finally we bound $\textbf{scan}_{m,n} (\bX)/mn - \bar{\bX}$ as a whole. By the definition of the scan statistic, 
	\beq
	\frac{\textbf{scan}_{m,n} (\bX)}{mn} - \bar{\bX} \geq \frac{Y_{S^*} (\bX)}{mn} - \bar{\bX}  =  \bar{\bX}_{S^*}  - \bar{\bX},
	\eeq
	here $\bar{\bX}_{S^*} = \sum_{(i,j) \in S^*} X_{ij} / mn$ represents the average in the submatrix indexed by $S^*$. Rewrite $X_{ij} = E(X_{ij}) + Z_{ij}$ for $X_{ij} \in S^*$, where $Z_{ij}$ has mean zero and bounded second moment. By Law of Large numbers, as well as (7) of \cite{arias2017distribution} (which is, $\E X_{ij} \geq \theta_\ddag$ for $X_{ij} \in S^*$),
	\beq
	\bar{\bX}_{\bS^*} = \frac{1}{mn} \sum_{X_{ij} \in S^*} \E X_{ij} + O_P(\frac{1}{\sqrt{mn}}) \geq \theta_\ddag + O_P(\frac{1}{\sqrt{mn}}).
	\eeq
	Note that $\bar{\bX} = O_P(1/\sqrt{MN})$ and $mn = o(MN)$, 
	\beq
	\frac{\textbf{scan}_{m,n} (\bX)}{mn} - \bar{\bX} \geq  \theta_\ddag + O_P(\frac{ 1}{ \sqrt{mn}} )
	\eeq
	By \eqref{condmain} we know that $\sqrt{mn} \theta_\ddag \rightarrow \infty$, or $\theta_\ddag \gg 1/\sqrt{mn}$, so we can rewrite the equation above as
	\begin{eqnarray} \label{eq:scan}
	\frac{\textbf{scan}_{m,n} (\bX)}{mn} - \bar{\bX} \geq  \theta_\ddag (1 + o_P(1)).
	\end{eqnarray}

	Plug in (\ref{eq:xmax}), (\ref{eq:var}) and (\ref{eq:scan}) into (\ref{eq:bern}), we have
	\beq \label{controlPValue}
	\mathfrak{P}(\bX)  \leq MN|\mathbb{S}_{mn}| \exp\bigg(- \frac{(1+o_P(1))\theta_\ddag^2mn}{2(1 + o_P(1)) + \frac{2}{c} (\log MN) \theta_\ddag(1+o_P(1))} \bigg).
	\eeq
	Assume in \eqref{condmain}, 
	\beq \label{highProbBound}
	\liminf \frac{\theta_\ddag^2 mn}{2 (m\log \frac{M}{m} + n\log \frac{N}{n})} \geq 1 + \epsilon
	\eeq
	with some constant $\epsilon > 0$. Then eventually with high probability we can bound the exponential part of \eqref{controlPValue} as 
	\beqn
	- \frac{(1+o_P(1))\theta_\ddag^2mn}{2(1 + o_P(1)) + \frac{2}{c} (\log MN) \theta_\ddag(1+o_P(1))} \leq -\frac{(1-\epsilon/2)\theta_\ddag^2mn}{2(1 + \epsilon/8)}.
	\eeqn
	Here we used the fact that the second term in the denominator in the exponent component is $o_P(1)$ from \eqref{extracond}. 
	Combined with \eqref{highProbBound}, we have the upper bound of $\log (\mathfrak{P}(\bX))$ as 
	\beqn
	\log (\mathfrak{P}(\bX)) &  \leq& \log (MN|\mathbb{S}_{mn}|) - \frac{(1-\epsilon/2)\theta_\ddag^2mn}{2(1 + \epsilon/8)}  \\ & \leq& \log ( MN|\mathbb{S}_{mn}| )  -(1+\frac{\epsilon}{4 + \epsilon/2}) (m\log \frac{M}{m} + n\log \frac{N}{n}) .
	\eeqn
	By the fact that $MN|\mathbb{S}_{mn} (X)| = (1+o(1)){N\choose n} {M \choose m}$ and $\log({N \choose n}) = (1 + o(1))n\log (N/n)$, we have eventually
	\beq
	\log (MN|\mathbb{S}_{mn}|) \leq ((1+\frac{\epsilon}{8 + \epsilon}))(m\log \frac{M}{m} + n\log \frac{N}{n}).
	\eeq
	So the log Bonferroni-corrected empirical p-value is bounded from above by
	\beq
	\log(\mathfrak{P}(\bX)) \leq -(\frac{\epsilon}{8 + \epsilon}) (m\log \frac{M}{m} + n\log \frac{N}{n}) 
	\eeq
	which finishes the proof.
	
	\subsubsection{Unidimensional permutation}
	
We refer to the proof of Theorem 2 in \cite{arias2017distribution}. Under the unidimensional permutation, from (17) and its following arguments of \cite{arias2017distribution}, for the permutation $p$-value with $(m,n)$ known, we have
\beq
\log \mathfrak{P}_{m,n} (\bX) \leq -(1+o_P(1)) \delta \log (|\mathbb{S}_{mn}|).
\eeq
Here $\delta$ is a positive constant that is less than the following quantity 
\beq
2\bigg[ \lim \inf \frac{\theta_\ddag \sqrt{mn}}{2(m\log(M/m) + n\log(N/n))} - 1 \bigg]
\eeq
Note that $\log (MN) = o(1)\log ( |\mathbb{S}_{mn}|)$ by \eqref{extracond}, therefore directly,
\begin{eqnarray}
\log \mathfrak{P} (\bX) &\leq& \log ( MN) + \log \mathfrak{P}_{m,n} (\bX) \nonumber \\ & \leq& -(1+o_P(1)) \delta \log  |\mathbb{S}_{mn}| + \log (MN) = -(1+o_P(1)) \delta \log ( |\mathbb{S}_{mn}|).
\end{eqnarray}

\subsection{Proof of \lemref{binary}}

	From \defref{kbinary}, 
\beq
	c' - c = \sum_{i = 0}^{\lfloor \log_2 c \rfloor - k} a_i(c) 2^i \leq \sum_{i = 0}^{\lfloor \log_2 c \rfloor - k}  2^i \leq 2^{\lfloor \log_2 c \rfloor - k + 1}.
\eeq 
	Observing that $c = 2^{\log_2 c}$, we have
\beq
	\frac{c - c'}{c} \leq  2^i \leq 2^{\lfloor \log_2 c \rfloor - k + 1 - \log_2 c } \leq 2^{1-k}.
\eeq

\subsection{Proof of \thmref{bonf2}}
We illustrate the bidimensional case here, since the unidimensional case is following the same proof strategy and basically a rework of the existing proof in \secref{proofbonf}. 

From \lemref{binary}, we may find $(m',n') \in S_{k_M}\times S_{k_N}$ such that $m' \leq m, n' \leq n$ and $m' = (1+o(1))m, n' = (1+o(1))n$. Now we consider performing permutation test on $(m',n')$ and bound $\mathfrak{P}_{m',n'}$ by the same way in \secref{proofbonf}. 

We rewrite (\ref{eq:bern}) as follows,
	\beq 
	  |S_{k_M}|| S_{k_N}||\mathbb{S}_{m'n'}| \exp\Bigg( - \frac{ m'n'(\textbf{scan}_{m',n'} (\bX)/m'n' - \bar{\bX})^2}{2\sigma^2_{\bX} + \frac{2}{3} (\bX_{\max} - \bar{\bX}) ( \textbf{scan}_{m',n'} (\bX)/m'n' - \bar{\bX})} \Bigg).
	 \eeq
Note that $|S_{k_M}|| S_{k_N}| \leq MN$. Combined with $m' = (1+o(1))m, n' = (1+o(1))n$, all we need to verify is a new version of (\ref{eq:scan}), which in details, we just need to show that 
\begin{eqnarray}  \label{bonf2scan}
\frac{\textbf{scan}_{m',n'} (\bX)}{m'n'} - \bar{\bX} \geq  \theta_\ddag (1 + o_P(1)).
\end{eqnarray}
This is done by realizing that $\textbf{scan}_{m',n'} (\bX) \geq Y_{S'}(\bX)$, where $S' \subset S^*$ and have $m'$ rows and $n'$ columns. Since $m',n' \to \infty$, the same argument yields 
\beq 
\frac{\textbf{scan}_{m',n'} (\bX)}{m'n'} \geq \theta_\ddag + O_P(\frac{1}{\sqrt{m'n'}}) \geq \theta_\ddag + O_P(\frac{1}{\sqrt{mn}}).
\eeq
Everything else follows the previous argument, and \eqref{bonf2scan} will be verified, which concludes the proof.

\subsection*{Acknowledgements}
The authors are grateful to Ery Arias-Castro for introducing this topic to our attention, and helpful discussions.

\bibliographystyle{chicago}
\bibliography{ref}

\end{document}